\magnification=1200
\parindent=0pt
\rightskip=0pt plus1fil
\input amssym.def
\input epsf
\scrollmode
\def \BZ {{\Bbb Z}}

\def \CK {{\cal K}}
\def \CC {{\cal C}}
\def \lcs {{\rm LCS}}
\def \ds {{\rm DS}}
\def \tag {\number \count1 .\number \count2 \global\advance\count2 by 1}
\def \bnum {\number \count3 \global\advance\count3 by 1}
\def \txt {\medskip (\tag )\ \ }
\def \rem {\medskip \noindent {\bf Remark \tag :} \ \ }
\def \notn {\medskip \noindent {\bf Notation \tag :}\ \ }
\def \theom {\medskip \noindent {\bf Theorem \tag :}\ \ }
\def \proof {\medskip \noindent {\it Proof:}\ \ }
\def \proofcont {\medskip}
\def \prop {\medskip \noindent {\bf Proposition \tag :}\ \ }
\def \cor {\medskip \noindent {\bf Corollary \tag :}\ \ }
\def \lemm {\medskip \noindent {\bf Lemma \tag :}\ \ }
\def \examp {\medskip \noindent {\bf Example \tag :}\ \ }
\def \qno {\leqno {(\tag )}}
\def \figr{\medskip \centerline {Figure \tag} \medskip}
\def\sqr#1#2{{\vcenter{\vbox{\hrule height.#2pt
     \hbox{\vrule width.#2pt height#1pt \kern#1pt
     \vrule width.#2pt} \hrule height.#2pt}}}}
\def\square{\ ${\mathchoice\sqr34\sqr34\sqr{2.1}3\sqr{1.5}3}$}
\def \cl#1{\overline {#1}}

\centerline {\bf VASSILIEV INVARIANTS AND KNOTS MODULO PURE BRAID SUBGROUPS}
\bigskip
\centerline {Theodore B. Stanford}
\centerline {Mathematics Department}
\centerline {United States Naval Academy}
\centerline {572 Holloway Road}
\centerline {Annapolis, MD\ \ 21402}
\medskip
\centerline {\tt stanford@nadn.navy.mil}
\bigskip

{\bf ABSTRACT}
\bigskip

We show that two knots have matching Vassiliev invariants of order less than $n$
if and only if they are equivalent modulo the $n$th group of the lower central
series of some pure braid group, thus characterizing Vassiliev's knot invariants
in terms of the structure of the braid groups.   We also prove some results
about knots modulo the $n$th derived subgroups of the pure braid groups, and
about knots modulo braid subgroups in general.

\count1 = 0
\count2 = 1
\bigskip
\bigskip
\noindent
{\bf
0.  INTRODUCTION
}
\bigskip

\txt
Vassiliev invariants [25] are a rather recent innovation in low-dimensional
topology.  The main reason that they have attracted a great deal of interest
is that they provide a unified framework in which to consider
``quantum invariants'' such as the Jones polynomial and its generalizations.
The main theorem of this paper
characterizes Vassiliev invariants in terms of the structure of the braid
groups, thus providing something of a
bridge between modern quantum
topology and classical topology and group theory:

\theom
If $K_1$ and $K_2$ are knots, then $v (K_1) = v(K_2)$ for
every Vassiliev invariant of order $<n$ if and only if there
exist a positive integer $k$ and braids $p,b \in B_k$ such
that $K_1 = \cl{b}$, $K_2 = \cl {pb}$,
and $p \in \lcs_n (P_k)$.

\txt
By $\lcs_n(P_k)$
we mean the $n$th group of the lower central series of the pure
braid group $P_k$, so $\lcs_n (P_k) = [\lcs_{n-1}(P_k), P_k]$
with $\lcs_1 (P_k) = P_k$.  By the closure $\cl{b}$
of a braid $b$ we mean the
standard notion where the $i$th strand on the bottom is looped
around and identified with the $i$th strand on the top.
Two good references on Vassiliev invariants, also
called finite-type invariants, are
Bar-Natan [1] and Birman [4].

\txt
For the purposes of this paper, we could define a Vassiliev
invariant as follows:  
Any braid closes up to form a unique link.  Therefore, any link
invariant pulls back to a braid invariant.
If the invariant takes
values in an abelian group 
(often an additive subgroup of the complex numbers)
then the invariant may be
extended to the group ring $\BZ B_k$ for all $k>0$.
Let $I_k$ be the augmentation ideal 
of $\BZ P_k$ (so $I_k$ is generated by all $p-1$ for $p \in P_k$).
A Vassiliev invariant of order $< n$ is then an invariant which
vanishes on $I_k^n$ for all $k$.
We formulate this more precisely and
relate it to the standard definition of Vassiliev invariant
in Proposition~2.17.

\txt
We define two knots $K_1$ and $K_2$ to be 
$\lcs_n$-equivalent
if there exist $k>0$ and braids $p,b \in B_k$ such that
$K_1 = \cl {b}$, $K_2 = \cl {pb}$,
and $p \in \lcs_n (P_k)$.  
We say they are $V_n$-equivalent if
$v(K_1) = v(K_2)$ for all Vassiliev invariants $v$ of order $< n$
taking values in any abelian group.
Then Theorem~0.2 is reworded in
Theorem~2.41 to state that two knots are 
$\lcs_n$-equivalent if and only
if they are $V_n$-equivalent.

\txt
Crucial inspiration for this paper was provided by Habiro's paper [11].
where a different answer is provided to the question of when two knots
have matching invariants of order $<n$.  
(In particular, the idea 
to set up the diagram of maps in Figure~2.2 was inspired
by Habiro.)
We do not, however, rely directly on any of Habiro's
results.  On the surface, his approach seems to have little to do
with ours.  For each posititive integer $n$, 
he defines a class of moves which involve twistings along
curves inside of handlebodies, and shows that these moves suffice
to go between any two knots with matching invariants up to order $n$.
However, Habiro also defines for each $n$ another move, which
he calls a ``$*^n$~move'', and this move may be shown to be equivalent
to the replacement (unoriented, see (0.9)) in a knot diagram
of $1 \in B_{n+1}$ by a particular $p \in \lcs_n (P_{n+1})$.
The braid $p$ may be written explicitly as a product of iterated
commutators, but it is easier to appeal to Theorem~4.2 of [19] to
show that it is in $\lcs_n (P_{n+1})$.

\txt
Another motivator and source of ideas for our work is Gusarov's
paper [10].  Gusarov was the first to construct a group structure
on knots dual to finite-type invariants. He also seems to be the first
to have realized the importance of invariants which are additive
under connected sum, and which are therefore homomorphisms from
his group into some other abelian group.
It is not hard to see that $\lcs_n$-equivalent knots in our sense
are $n$-equivalent in the sense of Gusarov and Ohyama [20],
and in particular that $\lcs_n$-trivial knots are $n$-trivial in their sense.
It follows from Theorem~0.2 that both notions of equivalence are
the same as $V_n$-equivalence.

\txt
The ``if" part of Theorem~0.2 was proved in [23], where we used it
to construct for any knot $K$ and positive integer $n$, an infinite
family of knots whose invariants matched those of $K$ up to order $n$.
There are two ways in which
Theorem~0.2 is actually a stronger converse to the result in [23]
than might have been hoped for.   First of all, to go between
two knots with the same invariants of order $<n$, one might expect
to have to make a number of replacements of $1 \in B_k$ with
$p \in \lcs_n (P_k)$ (not necessarily the same $k$ or 
$p$ each time), but
in fact what we show here is that the change can be done with
a single such replacement.  This is a consequence of
Proposition~1.22, that $\lcs_n$-equivalence is an equivalence
relation.  Note that we give here no way to control $k$.

\txt
In [23], a replacement of $1 \in B_k$ with $p \in \lcs_n(P_k)$ could
be made without regard for the orientation of strands.  That is,
the braid orientation did not need to coincide with the knot
orientation.  We showed that such a replacement did not change
invariants of order $<n$.  Theorem~0.2 says that if two knots have
equal invariants of order $<n$ then they are related by such
a replacement where the braid orientation does coincide with the
knot orientation.  This is the second way in which Theorem~0.2
is a strong converse to Theorem~1 of [23].

\txt
The lower central series has shown up 
elsewhere in the study of finite-type invariants.
Milnor defined a 
set of link invariants in [17]
using the lower central series of link groups.  Bar-Natan [2]
and Lin [14] showed that these invariants, properly interpreted,
are of finite type.  Kalfagianni and Lin [12] 
have given a characterization
of knots with trivial finite-type invariants of order $<n$ as
knots which admit a certain type of Seifert surface $F$ with certain curves
that push off into the $n$th group of the lower central series
of the group of the complement $F$ in $S^3$.   Garoufalidis and
Levine have shown [8] that the lower central series of the Torelli
group of a surface plays an important role in the study of
finite-type invariants of homology 3-spheres.

\txt
For groups $G$ whose abelianizations are cyclic, the lower central series stabilizes
after one term: $\lcs_n (G) = \lcs_2 (G)$ for $n>1$.  Therefore it is
not very interesting to study the lower central series of  knot groups or
braid groups (as opposed to link groups and pure braid groups).

\txt
A simple but important fact (Proposition~2.13) 
in the study of Vassiliev
invariants via braids is that the ideal of braids 
with $n$ double
points is the $n$th power of the augmentation ideal of $\BZ P_k$.
It is an almost trivial fact of group theory (Proposition~2.24)
that if $g$ is an element of any group $G$ and $g \in \lcs_n (G)$,
then $g - 1 \in I^n$, the $n$th power of the augmentation ideal
of $\BZ G$.
If the converse holds in a group $G$ then the group is said
to have the ``dimension subgroup property."  Not all groups
have this property (see Rips [21]), 
but pure braid groups do.  
Free groups do too, and in fact one way to show that
pure braid groups do is to use the result for free groups,
the semidirect product decomposition of $P_k$, and an argument of
Sandling [22].  See also Gupta [9]
for more discussion of the dimension subgroup problem.
For pure braids, therefore, being undetectable by invariants of order
$<n$ is exactly the same as being in $\lcs_n (P_k)$ (see also
Kohno [13]). 
Theorem~0.2 may then be viewed as extending this last statement
to a corresponding result about knots.
The proof of Theorem~0.2 in this paper does not use the 
dimension subgroup property of $P_k$, however
(nor do I know of a proof that does).

\txt
Falk and Randell [7] showed that $P_k$ is residually nilpotent.
That is,
$\cap_n \lcs_n (P_k) = \{1\}$ for all $k$,
and combining this with the dimension subgroup property of $P_k$
shows that Vassiliev invariants (of all orders) distinguish braids.
See also Bar-Natan [2].

\txt
Section~1 sets up some notation and basic results for a general
idea of $H$-equivalence of links, where $H = \{H_k\}$ is a sequence
of normal subgroups of $B_k$ satisfying an inclusion condition which
we call ``subcoherence'' (see Remark~1.6 for the reasoning behind this term).
The case of most interest to us is where
$H_k = \lcs_n (P_k)$.  We conclude the section
by showing (Theorem~1.41) that
knots modulo $\lcs_n$-equivalence form a group $\CK/\CK_n$
under connected sum.
We also show (Theorem~1.46) that for each $k$
there is a homomorphism of abelian groups from
$\lcs_n (P_k)/\lcs_{n+1}(P_k)$ to the group 
$\CK_n/\CK_{n+1}$ of $\lcs_n$-trivial knots
modulo $\lcs_{n+1}$-trivial knots.

\txt
The sequence
of results in Sections~1 and~2 is mostly self-contained.
For basic facts about groups, see Magnus, Karrass, and Solitar [15].
For basic facts about
braids and their closures, see Birman [3].
A slightly strengthened version of
Markov's Theorem on braid closures is used to prove Proposition~1.22,
that $\lcs_n$-equivalence is an equivalence relation.
A generalization of Alexander's theorem that any knot
may be written as a closed braid is
used in Section~2 to show
that the definition of Vassiliev invariant given
in (0.4) is the same as the standard one.

\txt
Section~2 is devoted to the proof of Theorem~0.2/2.41
We use the group of $\lcs_n$-equivalence classes from Section~1.
A key idea 
is to replace formal sums
of knots with connected sums.  
This is done with what we call composite relators, which
are of the form $K_1 \# K_2 - K_1 - K_2$, where $K_1$ and
$K_2$ are knots.   We show that the defining relations
for Vassiliev invariants, namely that any knot with $n$
double points (resolved in the standard way) vanishes, can
be generated by composite relators and relations of the
form $\cl {pb} = \cl {b}$ for $p \in \lcs_n (P_k)$,
the latter relations being precisely those that define $\lcs_n$-equivalence.
A consequence of this is Theorem~2.43, which says that if
two knots have matching additive invariants up to order $n$
then they have matching invariants up to order $n$.  This was
known for rational-valued invariants
(see Bar-Natan [1] or Gusarov [10]), since in this case
there is a Hopf algebra structure which shows that the space
of all Vassiliev invariants is the polynomial algebra generated
by the additive ones.  It is not clear whether Theorem~2.43 is
a nonvacuous extension, however, since all known invariants
taking values in $\BZ/q\BZ$ lift to integer-valued invariants.

\txt
Section~3 extends some of the results of Section~1 to the derived
series of $P_k$.  The $n$th group of the derived series is
given by $\ds_n (P_k) = [\ds_{n-1} (P_k), \ds_{n-1} (P_k)]$ 
with $\ds_1 (P_k) = P_k$.   We show that for all $n$, the
$\ds_n$-equivalence classes of knots form a group under connected sum,
and we show how to apply an argument from [23] to prove that 
every $\ds_n$-equivalence class of knots contains infinitely
many prime, alternating knots.

\bigskip
{\bf ACKNOWLEDGEMENTS}
\bigskip
I would like to thank J\'ozef Przytycki and Tatsuya Tsukamoto for 
introducing me to Habiro's work and for helpful discussions.
I would also like to thank Joan Birman, Kazuo Habiro, Jerry Levine,
Mark Kidwell, Hugh Morton, Rollie Trapp, and Arkady Vaintrob 
for helpful comments, suggestions,
and discusssions.

\count1 = 1
\count2 = 1
\bigskip
\bigskip
{\bf
1.  KNOTS MODULO BRAID SUBGROUPS
}
\bigskip

\txt
First we will develop some general ideas about $H$-equivalence of links.
Then we will consider the main case of interest to us, equivalence modulo
the lower central series of the pure braid groups $P_k$.
We will show that the equivalence classes of knots form a group $\CK / \CK_n$ under
connected sum, and define homomorphisms from the quotient groups 
of the lower central series of $P_k$ into $\CK / \CK_n$.

\notn
Let $B_k$ be the braid group on $k$ strands.  Let $\sigma_i$
be the standard generators of $B_k$, 
where $1 \le i < k$.  The defining relations
for $B_k$ are 
$$\sigma_i \sigma_j = \sigma_j \sigma_i \qquad \hbox {for $|i-j|>1$}
\qno$$
and 
$$\sigma_i \sigma_{i+1} \sigma_i 
= \sigma_{i+1} \sigma_i \sigma_{i+1} \qno$$
The symmetric group $S_k$ is the quotient of 
$B_k$ by the relations $\sigma_i^2 = 1$.  If $b$ is a braid then
its image under this quotient map to $S_k$ 
will be called the {\it permutation
associated} to $b$.
Let $P_k \subset B_k$ be the pure braid group.
A pure braid is one whose associated permutation is the identity.
For $b \in B_k$, denote 
by $\cl {b}$ the standard closure of $b$ into a knot or link.
For $k>0$, denote by $\iota_k$ the inclusion of 
$B_k$ into $B_{k+1}$
by adding an unbraided strand on the right end
(algebraicly, this means adding the generator $\sigma_k$).  We will sometimes 
identify a braid $b \in B_k$ with $\iota_k (b) \in B_{k+1}$,
or with a product of such $\iota_k$ into $B_m$ for $m>k$.

\notn
Let $H = \{H_k\}$ be a sequence of groups such 
that $H_k$
is a normal subgroup of
$B_k$ for each $k$, and such that 
$\iota_k (H_k) \subset H_{k+1}$ for each $k$.
We will call such an $H$ a 
{\it subcoherent sequence of braid subgroups}.
The subcoherent sequence $P_k \subset B_k$ will be denoted $P$.
If $H \subset P$ (that is, $H_k \subset P_k$ for all $k$),
then we will call $H$ a {\it subcoherent sequence of pure braid subgroups}.
If $H$ is a subcoherent sequence
of braid subgroups, then 
two links $L_1$ and $L_2$
will be called {\it $H$-equivalent} if there exist $k>0$, $b \in B_k$,
and $h \in H_k$ such that $L_1 = \cl {b}$ and $L_2 = \cl {hb}$.

\rem
Most of the sequences $H$ discussed in this paper, in particular
the groups of the lower central and derived series of $P_k$,
satisfy a stronger condition, which we might call {\it coherence}.  A coherent
sequence would satisfy $\iota_k (H_k) = H_{k+1} \cap \iota (B_k)$.  This
condition is equivalent to requiring that the map $B_k/H_k \to B_{k+1}/H_{k+1}$
induced by $\iota_k$ be injective.  Another reason one might want
to consider coherent sequences is that there is a
bijective correspondence between coherent sequences of subgroups and normal
subgroups of $B_\omega$, the direct limit of the $B_k$ under the~$\iota_k$.
Given a subcoherent sequence $\{H_k\}$, the direct limit $H_\omega$
is a normal subgroup of $B_\omega$.  Different subcoherent sequences may
have the same direct limit $H_\omega$, but for coherent sequences this limit is
unique.  Given the limit $H_\omega$, the coherent sequence may be recovered
by $H_k = B_k \cap H_\omega$.

\examp
It is not hard to see that two links are $P$-equivalent if and only
if they have the same number of components.

\txt
Birman and Wajnryb [6,26] treated Example~1.7 along with two other
interesting examples of subcoherent sequences.  
In fact, their sequences are coherent, and their
point of view is to consider the quotients $B_k/H_k$.

\notn
Given a group $G$ and two subsets $A,B \subset G$, let $[A,B]$ be the
subgroup of $G$ generated by all expressions 
$[a,b] = aba^{-1}b^{-1}$,
where $a \in A$ and $b \in B$.  If $H$ and $H^\prime$ are
sequences of braid subgroups,
then denote by $[H,H^\prime]$ the
sequence of subgroups
$[H_k, H^\prime_k] \subset B_k$.

\txt
The next proposition will allow us, in particular, to consider the sequences
$\lcs_n (P)$ and $\ds_n (P)$ (see Notations~1.40 and~3.8).

\prop
Let $H$ and $H^\prime$ be two subcoherent sequences of braid subgroups.
Then $[H,H^\prime]$ is a subcoherent sequence of braid subgroups.

\proof
Commutators are preserved by group homomorphisms.
\square

\examp
Consider $H = [P,P]$, so $H_k = [P_k, P_k]$
($= \lcs_2 (P_k) = \ds_2 (P_k) $).  
It turns out that all knots are $[P,P]$-equivalent, since by Theorem~0.2
this corresponds to having Vassiliev invariants of order $<2$
equal, and it is well-known that there are no Vassiliev knot invariants
of order $1$.   For links, $[P,P]$-equivalence will certainly detect
the number of components since $[P_k, P_k] \subset P_k$.
Moreover, it is not hard to see that for an $m$-component link
the unordered 
$m \choose 2$-tuple
of linking numbers is an
invariant of $[P,P]$-equivalence.
Note that this is a stronger invariant than the single $1$st-order
Vassiliev invariant of links,
which can be obtained by summing up all the pairwise linking numbers.
Thus Theorem~0.2 fails for links with more than two (undistinguished) components.

\txt
We will show below that $H$-equivalence is in fact an equivalence
relation, but first we need to describe some notation, and 
introduce a slight strengthening of the Markov Theorem from
Birman [3]. 
For the reader who is interested in keeping track of what gets used where, we
note that the sole purpose for the material from here to Proposition~1.22 is
to show that $H$-equivalence is an equivalence relation.  Without Proposition~1.22,
we could define $H$-equivalence to be {\it generated} by $\cl {b} = \cl {hb}$, with
$h$ and $b$ as in Notation~1.5.  This would weaken the result of Theorem~0.2
somewhat (see~(0.8)), but its proof would be essentially unchanged.

\txt
What we need is that a sequence of Markov moves
connecting two braids can be replaced by a sequence in which all the
stabilization is done first, followed by the destabilization.
Consider conjugation, the first Markov move.  Conjugate
braids in $B_k$ represent the same link under the closure operation,
and, moreover, the $b \in B_k$ in Notation~1.5 could be
replaced by any of its conjugates, since $h \in H_k$ may also
be replaced by any conjugate in $B_k$
without affecting the definition.
For notational convenience in what follows, 
we do not consider conjugation as a move itself,
but rather we consider the second Markov move,
stabilization, as a move between conjugacy
classes.  

\notn
Given a braid $b \in B_k$, denote
by $[b]$ its conjugacy class in $B_k$.  If $b_1 \in B_k$ and 
$b_2 \in B_{k+1}$, then we shall say that $[b_1]<[b_2]$ if
$[b_2]$ can be obtained from $[b_1]$ by a stabilization move.
This means that
there exist $b_1^\prime \in [b_1]$ and $b_2^\prime \in [b_2]$
such that $b_2^\prime = b_1^\prime \sigma_n$ or $b_2^\prime =
b_1^\prime \sigma_n^{-1}$.  This relationship induces a partial
ordering on the set of all conjugacy classes in all braid groups
$B_k$.  We will denote this partial order by $<<$. 

\txt
What the standard Markov Theorem says is that the equivalence relation
generated by $<$ is the same as equivalence of
the links obtained by closing up braids.  Two braid conjugacy classes
$[b]$ and $[c]$ will be equivalent if and only if there exists a
sequence $[b] = [b_0], [b_1], [b_2], \dots [b_p] = [c]$ 
such that for each $i$
either $[b_i] < [b_{i+1}]$ or $[b_{i+1}]<[b_i]$.  It is hard in
general to get any control over the sequence $[b_i]$.  The relation
$<<$ clearly generates the same equivalence relation as $<$, and
it turns out that with $<<$ it is never necessary to have more than
one intermediate conjugacy class between $b$ and $c$.  Both $<$ and
$<<$ behave very well with respect to $H$-equivalence, as shown in 
the next proposition. 

\prop  
Let $x \in B_k, y \in B_l$, and $h \in H_k$,
and suppose $[x] << [y]$.  Then there exists $j \in H_l$ such
that $[hx] << [jy]$.

\proof  It suffices to show the result for $[x] < [y]$,
where $x \in B_k$ and $y \in B_{k+1}$. Let
$b^{-1} yb = a^{-1} x a \sigma_k^\epsilon$, where
$\epsilon = \pm 1$.  Then $a^{-1} h x a \sigma_k^\epsilon
= a^{-1} h a b^{-1} y b = b^{-1} (b a^{-1} h a b^{-1}) y b$,
so we may set $j = b a^{-1} h a b^{-1}$. 
\square 

\txt
Here is the strengthened Markov Theorem:

\lemm 
Let $b_1 \in B_k$ and $b_2 \in B_l$ be such that 
$\cl {b_1} = \cl {b_2}$.  Then there exists a positive integer
$m$ and a braid $b_3 \in B_m$ such that $[b_3] >> [b_1]$ and
$[b_3]>>[b_2]$.

\proof  
It suffices, by the standard Markov theorem and induction, 
to show that if 
$b \in B_k$ and $c_1,c_2 \in B_{k+1}$ with $[c_1]>[b]$ and
$[c_2]>[b]$, then there exists
$d \in B_{k+2}$ such that $[d]>[c_1]$ and $[d]>[c_2]$.
Choose representatives of conjugacy classes so that 
$c_1 = b \sigma_k^{\epsilon_1}$ and 
$c_2 = \alpha^{-1} b \alpha \sigma_k^{\epsilon_2}$
where $\epsilon_1 = \pm 1$ and $\epsilon_2 = \pm 1$.
Now let $d = \alpha^{-1} \sigma_{k+1} \sigma_k^{\epsilon_1}
b \alpha \sigma_{k+1}^{-1} \sigma_k^{\epsilon_2}$.  
Using the braid relations (1.3) and~(1.4), we obtain
$$\eqalign {
\cl {\alpha^{-1} \sigma_{k+1} \sigma_k^{\epsilon_1}
b \alpha \sigma_{k+1}^{-1} \sigma_k^{\epsilon_2}}
&= \cl {\alpha^{-1} \sigma_{k+1} \sigma_k^{\epsilon_1} \sigma_{k+1}^{-1}
b \alpha \sigma_k^{\epsilon_2}} \cr
&= \cl {\alpha^{-1} \sigma_k^{-1} \sigma_{k+1}^{\epsilon_1} \sigma_k
b \alpha \sigma_k^{\epsilon_2}} \cr
&= \cl {\alpha^{-1} b \alpha \sigma_k^{\epsilon_2}}} \qno$$
and, conjugating by $\sigma_{k+1}$,
$$\eqalign {
\cl {\alpha^{-1} \sigma_{k+1} \sigma_k^{\epsilon_1}
b \alpha \sigma_{k+1}^{-1} \sigma_k^{\epsilon_2}}
&= \cl {\alpha^{-1} \sigma_k^{\epsilon_1} b \alpha
\sigma_{k+1}^{-1} \sigma_k^{\epsilon_2} \sigma_{k+1}} \cr
&= \cl {\alpha^{-1} \sigma_k^{\epsilon_1} b \alpha
\sigma_k \sigma_{k+1}^{\epsilon_2} \sigma_k^{-1}} \cr
&= \cl {\alpha^{-1} \sigma_k^{\epsilon_1} b \alpha} \cr
&= \cl {b \sigma_k^{\epsilon_1}}} \qno $$
\square

\prop  
Let $H$ be a subcoherent sequence of braid subgroups.
Then $H$-equivalence is an equivalence relation.

\proof  
The reflexive and symmetric 
properties are obvious.  Suppose that A, B, and C are links such
that $A$ is $H$-equivalent to $B$ and $B$ is $H$-equivalent to
$C$.  Then there exist $k$ and $l$, and braids $x,h \in B_k$,
and braids $y,j \in B_l$, with $h \in H_k$ and $j \in H_l$,
such that $A = \cl {xh}$, $B = \cl x = \cl y$, and $C =
\cl {jy}$.  By Lemma~1.19, there exists a braid $z \in B_m$
such that  $[z] >> [x]$ and $[z] >> [y]$.
By Proposition~1.17, there exist $h^\prime, j^\prime \in B_m$
such that $[h^\prime z] >> [hx]$ and $[j^\prime z] >> [jy]$.
Since $\cl {h^\prime z} = A$ and $\cl {j^\prime z} = C$, we get
that $A$ is $H$-equivalent to $C$.
\square

\txt
This paper is mostly concerned with pure braids and knots. 
If we want a notion of $H$-equivalence that concerns
only knots, then we need to assume that $H_k \subset P_k$.
Since
the closure of a pure braid isn't a knot, we choose a convenient
twist $t_k$ to put on the end of a pure braid so that it will close up
to a knot.  Once we have chosen such a twist, it is clear that any
knot can be represented as a pure braid plus the twist.
This twist will serve another useful purpose:
If $x$ is a $k$-strand braid which has been inserted into an
$m$-strand braid $y$, with $m>k$, and if we wish to
``shift $x$ to the right" in $y$, then the precise algebraic way to
do it is to replace $x$ by $t_m^{-1} x t_m$.

\notn  
Denote by $t_k \in B_k$ the braid 
$\sigma_{k-1}^{-1} \sigma_{k-2}^{-1} \dots \sigma_1^{-1}$.  
If $K_1$ and $K_2$ are knots, then we denote their
connected sum by $K_1 \# K_2$.

\prop
Let $x,y \in P_k$.  Then 
$$\cl {x t_k} \# \cl {y t_k} 
= \cl {x t_{2k}^{-k} y t_{2k}^{k+1}} \qno$$
where in the right-hand side of the equality $x$ and $y$ are identified
with their inclusions into $B_{2k}$.

\proof
Figure~1.27 shows the case $k=4$, which clearly generalizes.
\square

\epsfysize = 5truecm
\bigskip
\centerline {\epsffile {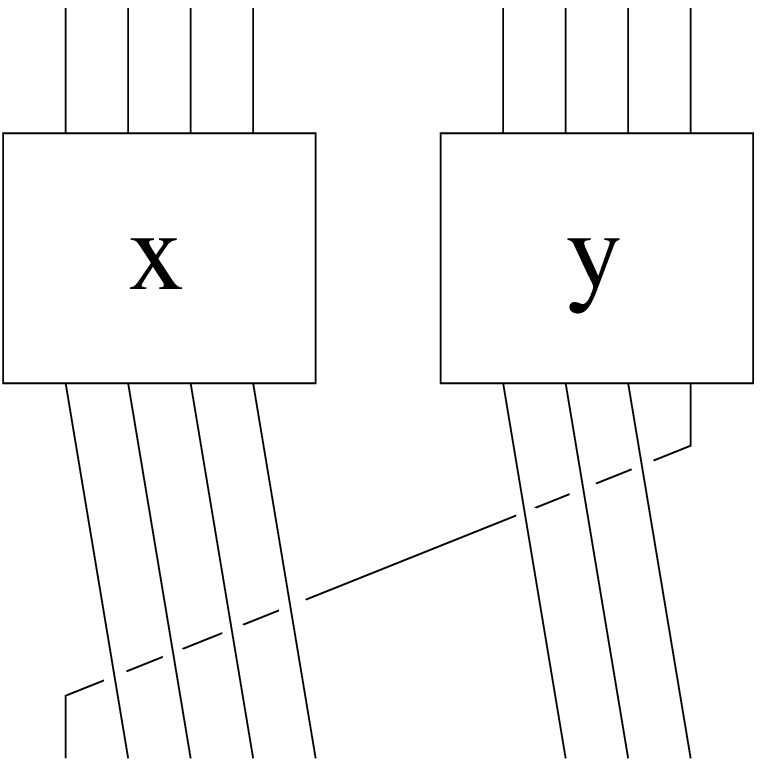} \hskip 1truecm
\vbox {\vskip1.8truecm \hbox {$xt_{2k}^{-k}yt_{2k}^k$} \vskip1.7truecm \hbox {$t_{2k}$} 
\vskip1truecm}}
\bigskip
\figr

\txt
The next two propositions say 
that $H$-equivalence is compatible with
the operation of connected sum.

\prop  
Let $K_1$, $K_1^\prime$, $K_2$, and $K_2^\prime$
be knots, and let $H$ be a subcoherent sequence of pure braid subgroups.
If $K_1$ is $H$-equivalent to $K_1^\prime$, and $K_2$ is
$H$-equivalent to $K_2^\prime$, then $K_1 \# K_2$ is $H$-equivalent
to $K_1^\prime \# K_2^\prime$.

\proof 
Without loss of generality we may assume that all four knots
are closures of braids $B_k$ for the same $k$.  Conjugating if
necessary, we may assume that the permutation 
associated to
each of these braids is the same as the permutation associated
to $t_k$.  So we may write 
$\cl {x_1 t_k} = K_1$, 
$\cl {h_1 x_1 t_k} = K_1^\prime$,
$\cl {x_2 t_k} = K_2$, and $\cl {h_2 x_2 t_k} = K_2^\prime$.  
By Proposition~1.25,
$$K_1 \# K_2 = \cl {x_1 t_{2k}^{-k} x_2 t_{2k}^{k+1}} \qno$$ 
and
$$K_1^\prime \# K_2^\prime 
= \cl {h_1 x_1 t_{2k}^{-k} h_2 x_2 t_{2k}^{k+1}}
= \cl  {h_1 t_{2k}^{-k} h_2 t_{2k}^k x_1 t_{2k}^{-k} x_2 t_{2k}^{k+1}} \qno $$
because $x_1$ commutes with $t_{2k}^{-k} h_2 t_{2k}^k$.
The result follows because $h_1 t_{2k}^{-k} h_2 t_{2k}^k \in H_{2k}$.
See Figure~1.27.
\square

\epsfysize = 7truecm
\bigskip
\centerline {\hbox {\epsffile {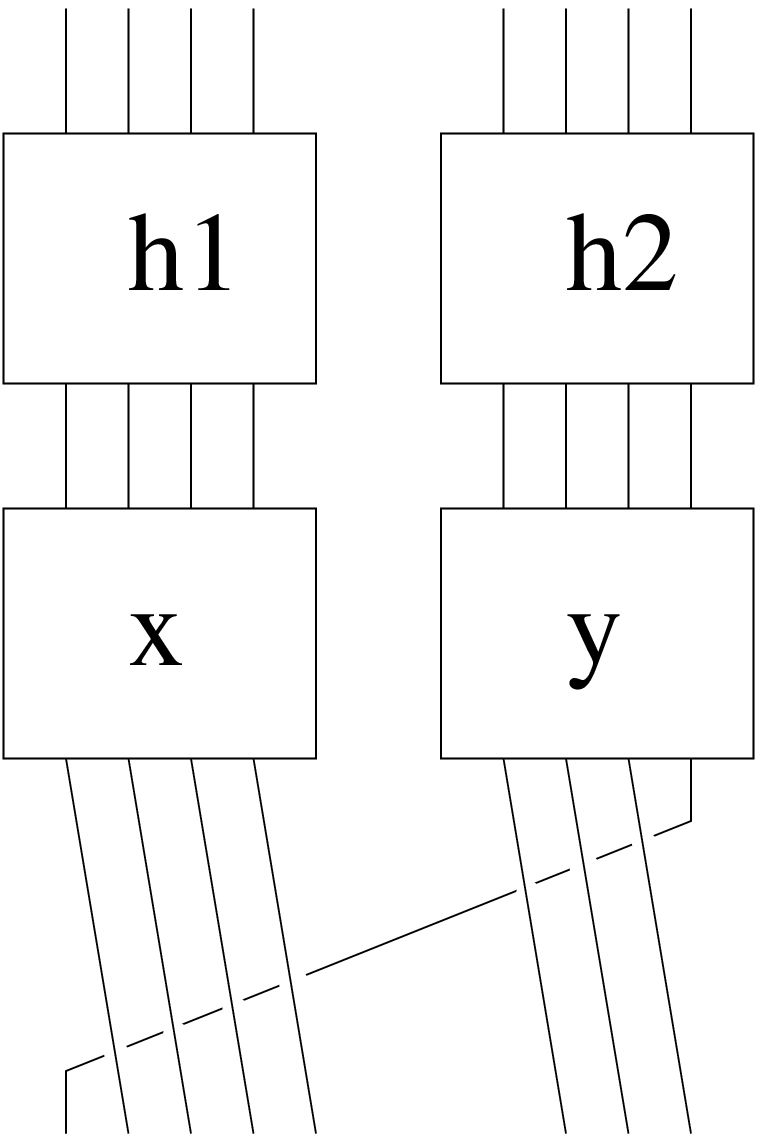} \hskip 1truecm
\vbox {\vskip2.4truecm \hbox {$h_1 t_{2k}^{-k} h_2 t_{2k}^k$} \vskip1.6truecm
\hbox {$xt_{2k}^{-k}yt_{2k}^k$} \vskip1.5truecm \hbox {$t_{2k}$} \vskip1truecm}}}
\bigskip
\figr

\prop  
Let $H$ be a subcoherent sequence of pure braid subgroups.
Then knots up to $H$-equivalence form an abelian monoid under the
operation of connected sum.

\proof  
Follows immediately from the fact that knots form an abelian
monoid, and Proposition~1.29.
\square

\txt
We are not only interested in a single sequence of subcoherent subgroups
$H$, but in various series of such sequences.  In this section and
the next we are mostly interested in the lower central series of
the pure braid group $P_k$.  In Section~3 we consider the derived series.

\txt
Proposition~1.37 allows us to relate the product of two braids with
the connected sum of two knots.  The proof is illustrated in
Figure~1.36. Conjugating the
first braid in the figure by $y$ (or, more precisely, by $t_8^{-3}yt_8^3$), we obtain
the second, which is then equivalent to the third modulo a commutator
involving $x$ and $y$.   

\epsfysize = 7truecm
\bigskip
\centerline {\epsffile {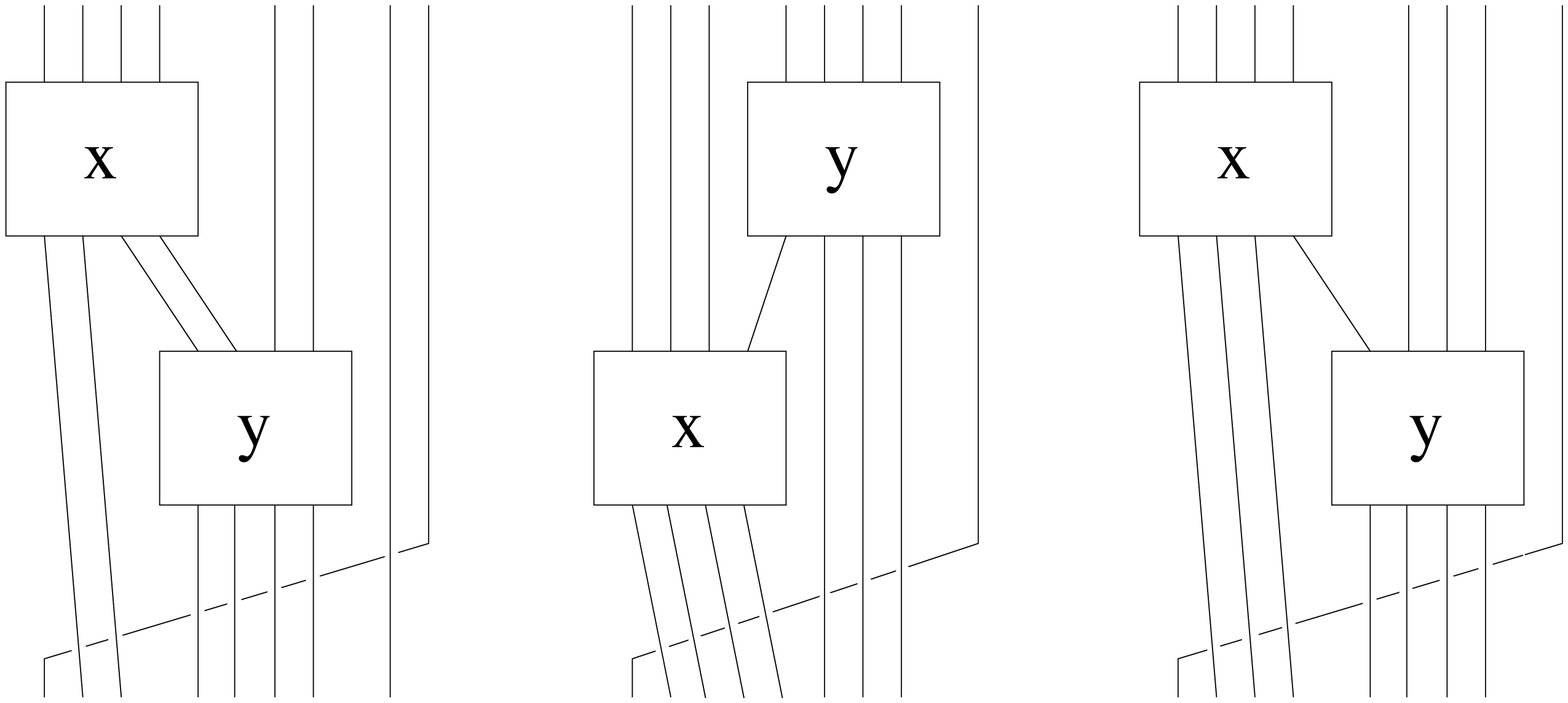}}
\bigskip
\figr

\prop
Let $H$ be a subcoherent sequence of pure braid subgroups.
Let $x \in H_k$ and $y \in P_k$.  
Then $\cl {x t_k} \# \cl {y t_k}$
is $[P,H]$-equivalent to $\cl {xy t_k}$.

\proof
Let $b_i = x t_{2k}^{-i} y t_{2k}^{i+1}$.
For $k=4$, the first braid in Figure~1.36 is $b_2$,
and the third is $b_3$.
When $i = 0$
then $\cl {b_i} = \cl {xy t_{2k}} = \cl {xy t_k}$, 
and when $i = k$ then
$\cl {b_i} = \cl {x t_k} \# \cl {y t_k}$ by Proposition~1.25.
We only need show that $\cl {b_i}$ is $[P,H]$-equivalent to 
$\cl {b_{i+1}}$.
Conjugating by $t_{2k}^{-i-1} y t_{2k}^{i+1}$,
we get
$$\cl {b_i} = \cl {x t_{2k}^{-i} y t_{2k}^{i+1}}
= \cl {t_{2k}^{-i-1} y t_{2k}^{i+1} x t_{2k}}
= \cl {[t_{2k}^{-i-1} y t_{2k}^{i+1},x]
x t_{2k}^{-i-1} y t_{2k}^{i+2}} \qno $$
which is $[P,H]$-equivalent to $\cl {b_{i+1}}$. 
\square

\prop 
If a knot $K$ is $H$-trivial, then there exists a knot
$K^\prime$ such that $K \# K^\prime$ is $[P,H]$-trivial.

\proof  
Let $K = \cl {h u t_k}$, where $h \in H_k$ and
$\cl {u t_k}$ is the unknot.
Set $K^\prime = \cl {h^{-1} t_k}$.
Then by Proposition~1.37, $K_1 \# K_1^\prime$ is $[P,H]$-equivalent
to $\cl {u t_k}$.
\square

\notn 
We denote by $\lcs_n (G)$ the $n$th group of the lower
central series of a group $G$.  
That is, $\lcs_1 (G) = G$, and $\lcs_n (G)
= [G,\lcs_{n-1}(G)]$ for $n>1$.  
If $H_k = \lcs_n (P_k)$, then we shall call $H$-equivalent knots
$\lcs_n$-equivalent and $H$-trivial knots $\lcs_n$-trivial.
Denote by $\CK$, $\CK_n$, and $\CK/\CK_n$ the abelian monoids of
all knots, all $\lcs_n$-trivial knots, and the $\lcs_n$-equivalence
classes of knots, respectively.

\theom  
For every $n>0$, the set $\CK/\CK_n$ of $\lcs_n$-equivalence classes
of knots forms a group
under connected sum.

\proof  
In light of Proposition~1.33, it suffices to show that for
every equivalence class modulo $\lcs_n$ there is an inverse
equivalence class. 
Let $K$
be a knot.  By Proposition~1.39, there exists a knot $K_1$ such that
$K \# K_1$ is $\lcs_2$-trivial.   Continuing inductively, we
obtain knots $K_i$ such that $K \# K_1 \dots K_{n-1}$ is $\lcs_n$-trivial.
Then $K_1 \# K_2 \dots \#K_{n-1}$ is the desired $\lcs_n$-inverse
for $K$.
\square

\txt
As remarked in (0.7), the group $\CK/\CK_n$ turns out to be the same as Gusarov's group
([10], see also [19]).

\rem
From a first glance at the proof of Theorem~1.41, it might appear
that the $\lcs_n$-equivalence class of a knot is related to the number
of prime factors it has.  This is not the case.  First of all, a knot
may be connected summed with many multiple copies of some knot in
$\CK_n$ without affecting its equivalence class.  Secondly, an argument
from [23] shows that any equivalence class can be represented by a prime
knot.  The argument is reworked in Section~3 for $\ds_n$-equivalence, which
implies the $\lcs_n$ case.

\notn
Denote by $\phi_k: P_k \to \CK$ the map given by 
$\phi_k (p) = \cl {p t_k}$. Denote also by $\phi_k$ any
quotient or restriction of this map.

\txt
The following Theorem is a straightforward consequence of Proposition~1.37.

\theom
For any $k,n>0$, The map $\phi_k : \lcs_n (P_k)/\lcs_{n+1} (P_k)
\to \CK_n / \CK_{n+1}$ is a homomorphism of abelian groups.

\rem
Note that the proposition implies that for any
$u \in P_k$ such that $\cl {u t_k}$ is the unknot,
$\phi (x) = \phi (xu)$.  Hence $\phi_k(x)$ may be evaluated by
closing up $xv$ (where $v = ut_k$) for any $v \in B_k$ such that $\cl v$
is the unknot and $v$ has the same associated permutation as $t_k$.

\txt
The following Theorem is an immediate consequence of Theorem~1.11 in [19].
It says that all $\lcs_n$-equivalence classes are realized by braids
of index $\le n$.

\theom
For any $n>0$, the homomorphism 
$\phi_{n+1}: \lcs_n (P_{n+1})/\lcs_{n+1} (P_{n+1}) \to \CK_n / \CK_{n+1}$
is surjective.

\rem
A good understanding of the kernel of 
$\phi_k : \lcs_n(P_k)/\lcs_{n+1} (P_k) \to \CK_n/\CK_{n+1}$
could lead to an alternate proof of Kontsevich's Theorem that
any weight system can be integrated to give a Vassiliev invariant.
See Bar-Natan [1].  Kontsevich's Theorem applies only to rational-valued
invariants, and perhaps this could be generalized to torsion invariants.

\rem
Since $\phi_{k+1} \circ \iota_k = \phi_k$, it is possible to define
a unique $\phi_\omega: P_\omega \to \CK$, 
where $P_\omega$ is the direct limit of the $P_k$ under the
inclusions $\iota_k$ as in Remark~1.6.
For any $k$,
restricting $\phi_\omega$ to $P_k$ will yield $\phi_k$.

\count1 = 2
\count2 = 1
\bigskip
\bigskip
{\bf
2. VASSILIEV INVARIANTS
}
\bigskip

\txt
The purpose of this section is to prove that equivalence of knots
modulo $\lcs_n(P)$ is the same as equivalence modulo Vassiliev
invariants of order $<n$.  The proof is summarized in the following
commutative diagram, which we describe briefly before we begin with
the details.

\bigskip
\centerline {\vbox {
$$

\def\mapright#1{\smash{\mathop{\longrightarrow}\limits^{#1}}}
\def\mapdown#1{\Big\downarrow\rlap{$\vcenter{\hbox{$\scriptstyle#1$}}$}}
\matrix
{\BZ \CK & \mapright{\alpha_n} & \BZ(\CK/\CK_n) & \mapright {\beta_n} & \CK/\CK_n \cr
\Big\Vert & & \mapdown{} & & \mapdown \cong \cr
\BZ \CK & \mapright{a_n} & V_n & \mapright {b_n} & W_n \cr}
$$}}
\figr

\txt
All maps in this diagram are considered as $\BZ$-module homomorphisms
(though some have additional structure).  
All maps 
except the vertical one in the middle are surjective.
The $\BZ$-modules are finitely-generated, except for $\BZ \CK$ and
$\BZ (\CK/\CK_n)$.
$\BZ \CK$ is the $\BZ$-module
freely generated by all knot types.  
The group $\CK/\CK_n$ is the 
group of $\lcs_n$-equivalence classes of
knots, as in Section~1, and $\BZ (\CK/\CK_n)$
is its integral group ring.
The map $\alpha_n$ is the map
which associates to each knot its $\lcs_n$-equivalence class, 
extended linearly to $\BZ \CK$.
The map $a_n$ is the quotient
of $\BZ \CK$ by all knots with $n$ double points, 
considered as an alternating sum of $2^n$ nonsingular knots in the
usual way.  Therefore, a Vassiliev knot invariant is a homomorphism
from $V_n$ to some abelian group.  
The map $b_n$ is the quotient of $V_n$ by all relations
$K_1 \# K_2 = K_1 + K_2$, where the $K_i$ are the images of 
knots (not singular knots or linear combinations of knots)
under $a_n$.  An additive Vassiliev invariant $w$, one in
which $w (K_1 \# K_2) = w (K_1) + w(K_2)$, is just a homomorphism
from $W_n$ to some abelian group.  
The map $\beta_n$ 
sends $K \in \CK/\CK_n$ to $K \in \CK/\CK_n$ and is extended linearly to
$\BZ (\CK/\CK_n)$.
The main theorem is proved by
showing that the kernel of $\beta_n \circ \alpha_n$ is equal to
the kernel of $b_n \circ a_n$, and hence we get
the isomorphism shown between
$\CK/\CK_n$ and $W_n$.

\txt
The vertical arrow in the middle will not be given a name.  It is there to
illustrate Proposition~2.29, which states that $a_n$ factors through $\alpha_n$,
or, in other words, that $\lcs_n$-equivalent knots have equal Vassiliev invariants
of order $<n$.  

\rem
We note that if $GK$ is defined to be the Grothendieck group of
knots under connected sum, then $GK$ fits into the diagram above
with $\BZ \CK$ mapping to $GK$ and $GK$ mapping to $\CK/\CK_n$ and
$W_n$.  The point is that the relations
that define the maps $b_n$ and $\beta_n$
are the images under $a_n$ and $\alpha_n$ of the relations
that take $\BZ \CK$ to $GK$.

\txt
Now we will develop the modules and maps in Figure~2.2 more formally.

\notn   
We will denote the set of all smooth, tame, and oriented
knots in oriented $S^3$ by $\CK$.  Since this set is an abelian monoid under
connected sum $\#$, we will sometimes use $1 \in \CK$ to refer to the unknot.
The $\BZ$-module freely generated by
$\CK$ will be denoted $\BZ \CK$.  The connected sum operation $\#$ extends
linearly to $\BZ \CK$.
A crossing marked with a dot
or vertex in a knot diagram will be called a {\it double point}, and
will be regarded here as simply
a shorthand notation for the difference (in $\BZ \CK$) 
between the knot  with the positive crossing and the knot
with the negative crossing.  So $K_\times = K_+ - K_-$, where
the local diagrams are shown in Figure~2.8.
A double point in a braid diagram will serve the same purpose, being
the difference in $\BZ B_k$ between two braids which are the same
except that one has a positive crossing at the double point and
the other has a negative crossing.
We will denote
by $\CC_n$ the subgroup of $\BZ \CK$ generated by all knot
diagrams with $n$
double points.  We will denote by $a_n: \BZ \CK \to V_n$
the quotient map with kernel $\CC_n$.  
Let $\alpha_n$ be the map from $\CK$ to $\CK/\CK_n$
which sends each knot to its $\lcs_n$ equivalence class.
Denote also by $\alpha_n$ the $\BZ$-linear extension
to $\BZ \CK \to \BZ (\CK/\CK_n)$

\epsfxsize = 12truecm
\bigskip \centerline {\epsffile {fig4.eps}}
\medskip \hskip 2.8truecm $K_\times$ \hskip 4.4truecm $K_+$ \hskip 4.4truecm $K_-$
\figr

\prop  
Let $A$ be an abelian group (written multiplicatively)
with the usual $\BZ$-module structure (given by $z(a) = a^z$
for $a \in A$ and $z \in \BZ$),
and let $\BZ A$
be its integral group ring with the usual $\BZ$-module structure 
(given by
$z \sum z_i a_i = \sum (zz_i) a_i$).
Let $\beta: \BZ A \to A$ be defined by
$\beta (a) = a$ and extended to $\BZ A$ linearly.
Then the kernel of $\beta$ is generated by
expressions of the
form $a_1 a_2 - a_1 - a_2$ for $a_1, a_2 \in A$.

\proof
Note that $\BZ A$ is freely generated by $a \in A$
as a $\BZ$-module, so $\beta$ is a well-defined
$\BZ$-linear map. We have $\beta (a_1 a_2 - a_1 - a_2)
= a_1 a_2 a_1^{-1} a_2^{-1} = 1_A$ for $a_i \in A$.  
Conversely, suppose that 
$\beta (\sum z_i a_i) = \Pi a_i^{z_i} = 1_A$.
We may reduce $\sum |z_i|$ to $1$ or $0$
by applying $a_1 + a_2 = a_1 a_2$
and $a_1 - a_2 = a_1 a_2^{-1}$ repeatedly.  Then if there is one
remaining element
of $A$ in the sum it must be $1_A$.
\square

\notn
An expression of the form 
$K_1 \# K_2 - K_1 - K_2$,
where $K_1$ and $K_2$ are knots, will be called a {\it composite relator}.
We will also refer to the image of such an expression under
$a_n: \BZ \CK \to V_n$ or 
$\alpha_n: \BZ \CK \to \BZ (\CK/\CK_n)$ as a composite relator.
We will denote by $b_n : V_n \to W_n$ the quotient map whose kernel
is generated by all composite relators in $V_n$.
We will denote by $\beta_n: \BZ (\CK/\CK_n) \to \CK/\CK_n$ the $\BZ$-module
map which is the identity when restricted to the subset $\CK/\CK_n$ of
$\BZ \CK/\CK_n$.

\prop
If $x,y \in \CC_1 \subset \CK$, then $\beta_n \circ \alpha_n (x \# y)
= b_n \circ a_n (x \# y) = 0$.

\proof
First note that the unknot is a composite relator (let $K_1 = K_2$ be the unknot
in the definition).
Both $x$ and $y$ are $\BZ$-sums of
expressions of the form $K-1$, where $K$ is
a knot.  So the Proposition~reduces to the fact that $(K_1-1)\#(K_2-1)$ is
the sum of the composite relators $1$ and $K_1\#K_2 - K_1 - K_2$.
\square

\notn  
Let $G$ be a group, and $\BZ G$ its integral group
ring.  We will denote by $I$ the augmentation ideal of $\BZ G$,
which is the kernel of the map to $\BZ$ which sends each group
element to $1$.  $I$ is generated (as a $\BZ$-submodule) by
all $g-1$ where $g \in G$.  We will denote by $I_k$ the augmentation
ideal of $P_k$.

\prop
$I_k^n$ is generated by all $k$-strand braid diagrams with $n$ double points.

\proof
A diagram with $n$ double points is of the following form in $\BZ P_k$:
$$w_1 (\sigma_{i_1}-\sigma_{i_1}^{-1}) w_2 (\sigma_{i_2}-\sigma_{i_2}^{-1})
\dots w_n (\sigma_{i_n}-\sigma_{i_n}^{-1}) w_{n+1} \qno $$
where $w_j \in B_k$.  Set $v_j = w_1 \sigma_{i_1}^{-1} w_2 \sigma_{i_2}^{-1} 
\dots w_n \sigma_{i_n}^{-1}$.  Then (2.14) may be rewritten as
$$(v_1 \sigma_{i_1}^2 v_1^{-1} -1) (v_2 \sigma_{i_2}^2 v_2^{-1} - 1) \dots
(v_n \sigma_{i_n}^2 v_n^{-1} - 1) v_n w_{n+1} \qno $$
Since $v_j \sigma_{i_j}^2 v_j^{-1} \in P_k$, (2.15) is an element of $I_k^n$.

\proofcont
Since $I_k$ is generated as a $\BZ$-module by elements of the form $x-1$,
it suffices to show that an element of the form $(x_1-1) (x_2-1) \dots (x_n-1)$,
where $x_i \in P_k$, is a linear combination of diagrams with $n$ double points.
Then, since each $x_i$ is a pure braid, each may be undone
by crossing changes (or, to put it more algebraicly, $P_k$ is normally generated
by $\sigma_k^2$).
This gives each $x_i -1$ as the $\BZ$-sum of
diagrams with one double point each.  Putting them all together with
the distributive property gives the desired result. 
\square

\notn  
Let $t_k \in B_k$ be the braid
$\sigma_{n-1}^{-1} \sigma_{n-2}^{-1} \dots \sigma_1^{-1}$
as in the Section~1.
As before, we will denote by $\cl x$ the standard closure of
a braid into a knot or link.  For each $k$ 
we may extend the closure operation
$\BZ$-linearly to a map from $\BZ B_k$ to $\BZ \CK$.

\prop  
$\CC_n$ is generated by all elements of the form
$\cl {x t_k}$, where $k>0$ and $x \in I_k^n$.

\proof  
Any knot diagram with $n$ double points may be written as
a closed braid diagram with $n$ double points.  This is
a generalization of Alexander's Theorem that every link
can be written a closed braid, and is proved in Birman [4].
After conjugating
by an element of $B_k$
if necessary, we may assume that the permutation associated to
$c$ is the same as that associated to $t_k$, and we may write $c = q t_k$,
where $q$ is a pure braid diagram with $n$ double points.  Proposition~2.13
completes the proof.
\square

\notn
A {\it relator of length $m$ and order $\ge n$} is
an expression of the form
$$\cl {(x_1-1)(x_2-1) \dots (x_m - 1) y t_k} \qno $$ 
where $k>0$ and
$x_i \in \lcs_{n_i} (P_k)$ for all $i$, and $n = n_1 + n_2 \dots +n_m$.

\txt
Since $\cap_n \lcs_n (P_k) = \{1\} \subset P_k$, it would in fact be possible
to define the exact order of a relator by requiring $x_i \notin \lcs_{n_i+1} (P_k)$,
but this is not necessary for our purposes.  These relators are generalizations
both of the generators of $\CC_n$ and of $\lcs_n$-equivalence, and will be our
tool for moving between the two.

\prop 
$\CC_n$ is generated by all relators of length $n$
and order $\ge n$.

\proof
This is just a restatement of Proposition~2.17, using Notation~2.18.
\square

\prop  
The kernel of $\alpha_n$ is generated by relators
of length 1 and order $\ge n$. 

\proof  
This is immediate from the definition of $\lcs_n$-equivalence.
\square

\txt
The next two propositions are well-known facts of group theory.

\prop
Let $G$ be a group and let $I \subset \BZ G$ be its augmentation ideal.
If $x \in \lcs_n (G)$, then $x - 1 \in I^n$.

\proof
Let $x \in \lcs_n(G)$ and $y \in G$.  Assume inductively that therefore $x - 1 \in I^n$.
We have 
$[x,y]-1 = (xy - yx) x^{-1}y^{-1} = ((x-1)(y-1) - (y-1)(x-1)) x^{-1}y^{-1} \in I^{n+1}$.
Also, if $x-1 \in I^n$ and $y-1 \in I^n$, then $xy-1 = x (y-1) + (x-1) \in I^n$;
and if $x-1 \in I^n$ then $x^{-1} - 1 = -x^{-1} (x-1) \in I^n$.
\square

\prop
Let $G$ be a group.  Let $x \in \lcs_m (G)$ and $y \in \lcs_n (G)$.  Then 
$[x,y] \in \lcs_{m+n} (G)$.

\proof
It is not hard to find commutator identities, which hold in a free group
and hence in any group, and which allow one to break up and reorder
commutators of products and commutators of commutators.  For example,
$$[xy,z] = (x[y,z]x^{-1})[x,z] \qno$$
and
$$[[x,y],z] = (xyx^{-1}z[[z^{-1},y^{-1}],x]z^{-1}xy^{-1}x^{-1})(x[y,[x^{-1},z]])x^{-1}) \qno$$
\square

\txt
The next proposition was proved (with different notation) in [23].

\prop
The map $a_n$ factors through the map $\alpha_n$.

\proof
By Propositions~2.21 and~2.22, this is the same as saying that a relator of
length 1 and order $n$ can be written as a 
sum of relators of length
$n$ and order $n$.  This follows from Proposition~2.24.
\square

\prop  
A relator of length $m$ and order $\ge mn$ is
a linear combination of relators of length $1$ and order $n$.

\proof
In such a relator,
there must exist $i$ such that $x_i \in \lcs_n (P_k)$.  Using the
distributive property of multiplication, we may write our given
relator as a sum of expressions of the form 
$\cl {w (x_i -1) y t_k} = 
\cl {(w x_i w^{-1} - 1) w y t_k}$, where
$w,y \in P_k$.
\square

\txt
The following proposition is the key to the theorems of this section.
Its proof depends on the same idea as Proposition~1.37.
There a connected sum was made into a braid product by sliding
a braid $y$ around and around modulo the next term down in the lower
central series.  Here we do the opposite.  Given a product of
expressions $(x_i-1)$ in $\BZ P_k$, we can slide $x_1 -1$ around
and around, modulo relations closer to what we want, until
it forms a connected sum which is zero because of the composite
relators.

\prop  
Any relator of order $\ge n$ can be written as the
$\BZ$-sum of
composite relators and 
relators of length $1$ and order $\ge n$.

\proof  Suppose on the contrary that there is some relator of order
$\ge n$ which cannot be so written.  Choose such a relator with
minimal length $m$, and among those with minimal length, one with
maximal order.  The order must be less than $mn$ by Proposition~2.30.
Let us write our relator 
$\cl {(x_1-1)(x_2-1) \dots (x_m - 1) y t_{2k}}$, where initially
$y \in P_k$ and
$x_i \in \lcs_{n_i} (P_k)$, included into $B_{2k}$ in the standard way.
As the argument progresses, the $x_i$ and $y$ will shift and and combine
to occupy more than $k$ strands in $P_{2k}$.

\proofcont
First of all, we may interchange $(x_i -1)$ with $(x_{i+1}-1)$
modulo relators of shorter length or greater order:
$$\eqalign{
(x_i - 1)(x_{i+1} - 1) -   (x_{i+1} - 1)(x_i - 1) 
&= x_i x_{i+1} - x_{i+1} x_i \cr 
&= ([x_i, x_{i+1}] - 1) + ([x_i, x_{i+1}] - 1)(x_{i+1} x_i - 1)} \qno $$
where the relation containing $([x_i, x_{i+1}] - 1)$
has the same order
but shorter length than the original, and the relation containing
$([x_i, x_{i+1}] - 1)(x_{i+1} x_i - 1)$ has the
same length but greater order.  We are applying Proposition~2.25 here.
Secondly, 
$x_m-1$ may be exchanged with
$y$ modulo relators of greater order:
$$(x_m - 1)y - y (x_m - 1) = ([x_m, y] - 1)yx_m \qno$$
Thirdly, by conjugation we have
$$\eqalign {
& \cl {(x_1-1)(x_2-1) \dots y (x_m - 1) t_{2k}} \cr
&= \cl {(t_{2k}^{-1}x_mt_{2k} - 1)(x_1-1)(x_2-1) \dots
(x_{m-1}-1) y t_{2k}}} \qno$$
Finally, applying (2.33), (2.34), and (2.35) repeatedly, we obtain
$$\eqalign {
& \cl {(x_1-1)(x_2-1) \dots (x_m - 1) yt_{2k}} \cr
&= \cl {t_{2k}^{-k}(x_1 - 1)t_{2k}^k
(x_2-1)(x_3-1) \dots (x_{m-1}-1) y t_{2k}}} \qno$$
modulo relators of shorter length and relators of
greater order. This last
relator is equal to 
$$\cl {(x_2-1)(x_3-1)\dots (x_m-1) y t_k} \# \cl {(x_1 - 1) t_k} \qno$$
by Proposition~1.25, and
is therefore a $\BZ$-sum of composite relators by Proposition~2.11.
\square

\prop  
The kernel of $\beta_n \circ \alpha_n: \BZ \CK \to \CK/\CK_n$ is equal
to the kernel of $b_n \circ a_n: \BZ \CK \to W_n$.

\proof  
By the preceding propositions, both kernels are generated
by relators of length 1 and order $\ge n$ and by composite
relators.
\square

\notn
We say that two knots are $W_n$-equivalent if they are equivalent in $W_n$.
We say that they are $V_n$-equivalent if they are equivalent in $V_n$.

\prop  
If two knots are $W_n$-equivalent, then they
are $\lcs_n$-equivalent.

\proof
In both cases, this just says that the difference of the two knots is
in the kernel of $b_n \circ a_n$.
\square

\theom
Two knots are $V_n$-equivalent if and only if they are $\lcs_n$-equivalent.

\proof
If two knots are equivalent modulo $\CK_n$ then they are $V_n$-equivalent
by Proposition~2.29.
If two knots are $V_n$-equivalent  then they are $W_n$-equivalent,
and therefore they are $\lcs_n$-equivalent.  
\square

\txt
As noted in (0.16), it is not known whether the following theorem
actually applies to anything new.

\theom
If $K_1$ and $K_2$ are knots, then $v(K_1) = v(K_2)$ for all Vassiliev
invariants of order $<n$ if and only if $w (K_1) = w(K_2)$ for all
additive Vassiliev invariants of order $<n$.

\count1 = 3
\count2 = 1
\bigskip
\bigskip
{\bf
3.  THE DERIVED SERIES
}
\bigskip

\txt
If instead of the set $\CK$ of all knots we restrict ourselves
to $H$-trivial knots for some $H$, we get the following
immediate generalization of Theorem~1.41:

\theom
Let $H = \{H_k\}$ be a subcoherent sequence of pure braid subgroups.
Let $H_k^{(n)} = [P_k, H_k^{(n-1)}]$, with $H_k^{(1)} = H_k$.
Then for any $n>0$, the set of $H^{(n)}$ equivalence classes
of $H$-trivial knots forms a group under connected sum.

\txt
Given $H$, we would like to consider $H$-trivial knots modulo
$\lcs_n (H)$, but what stops us is that we have no control over
the braid $u$ in the proof of Proposition~1.39.  So we will
consider a stronger notion of $H$-triviality, somewhat artificial
but useful nonetheless.

\notn
Let $H$ be a subcoherent sequence of pure braid subgroups.
We shall say that a knot $K$ is {\it strongly $H$-trivial} if
there exist $k>0$ and $h \in H_k$ such that 
$K = \cl {h t_k}$.

\prop
Let $H$ be a subcoherent sequence of pure braid subgroups.
Then the connected sum of two strongly $H$-trivial knots is strongly
$H$-trivial.

\proof
Use Proposition~1.25, with $x,y \in H_k$.
\square

\txt
Recall that Proposition~1.11 tells us that we can form the commutator
between two subcoherent sequences of braid subgroups and obtain a
new subcoherent sequence.

\prop
Let $K$ be a strongly $H$-trivial knot.  Then there exists a strongly
$H$-trivial knot $K^\prime$ such that $K \# K^\prime$ is strongly
$[H,H]$-trivial.   

\proof
Follow the proofs of Propositions~1.37 and~1.39.  The difference
is that there is no $u$ in this case, so that $x$ and $y$ may
be assumed by induction to be in $\ds_n (P_k)$.
\square

\notn
The $n$th group of the derived series of a group $G$ will be denoted $\ds_n (G)$.
Set $\ds_1 (G) = G$, and $\ds_{n+1} (G) = [\ds_n (G), \ds_n (G)]$.  
If $H = \{H_k\}$ is a sequence of groups, then $\ds_n(H)$ will denote
the sequence $\{\ds_n(H_k)\}$. If two knots
are $\ds_n(P)$-equivalent, then we shall call them {\it $\ds_n$-equivalent}, and
a $\ds_n(P)$-trivial knot will be called {\it $\ds_n$-trivial}.

\theom
Let $H$ be a subcoherent sequence of pure braid subgroups.
Then for any $n>0$, the set of $\ds_n (H)$-equivalence classes of
strongly $H$-trivial knots forms a group under connected sum.

\proof
Same as Theorem~1.41.
\square

\cor
For any $n>0$, the $\ds_n$-equivalence classes of knots form a group
under connected sum.

\proof
If $H_k = P_k$, then all knots are strongly $H$-trivial.
\square

\txt
Corollary~3.10 implies Theorem~1.41, since $\ds_n (G) \subset \lcs_n (G)$
for any $n$ and any group $G$.  More generally, we obtain

\cor
Let $H$ be a subcoherent sequence of pure braid subgroups.
Then for any $n>0$, the set of $\lcs_n (H)$-equivalence classes of strongly
$H$-trivial knots forms a group under connected sum.

\rem
It follows from Proposition~2.25
that for any group $G$, $\ds_n (G) \subset \lcs_{2^{n-1}} (G)$.
At the very least, then, the $\ds_n$-equivalence classes of knots
pick up Vassiliev invariants much faster than the $\lcs_n$-equivalence
classes.  It would be interesting to know whether they pick up
anything else.  Note that the quotients $\lcs_n (P_k)/\lcs_{n+1} (P_k)$
are all finitely-generated abelian groups, wherease the quotients
$\ds_n (P_k)/\ds_{n+1} (P_k)$ are abelian but are not finitely-generated.
So on the one hand it seems very possible that $\ds_n$-equivalence
is more than just a faster way to get at Vassiliev invariants, and on
the other hand it seems that it will be very difficult to analyze
$\ds_n$-equivalence. 

\txt
In [23] we showed that for any $n>0$ and any $\lcs_n$-equivalence class
of knots, it is possible to choose a representative of that class
which is prime and alternating.  In fact, it is possible to choose
an infinite number of such representatives, and moreover, the
same trick works for $\ds_n$-equivalence classes, as we shall
show in Theorem~3.23 below.  
The idea for the following approach to $B_3$ came from Birman and Menasco [5].

\notn
In $B_3$, let $a = \sigma_1$ and $b = \sigma_2$.  Let $A = a^{-1}$
and $B = b^{-1}$.  Let $d = bab = aba$ and $D = ABA = BAB$.
Let $w$ stand for a word in $a$ and $B$ (no $A$ or $b$), not
necessarily the same word each time it is used.  Let $W$
stand for a word in $A$ and $b$, not necessarily the same each
time. 

\prop
$DaD = dad = b$, $DbD = dbd = a$, $DAD = dAd = B$, $dBd = DBD = A$,
$DWD = dWd = w$, $DwD = dwd = W$.

\proof
Follows from the braid relations (1.4).
\square

\prop
For any $n>0$, there exist in $\ds_n (P_3)$ words of each of
the following forms:  $awa, awB, Bwa, BwB, awaD, awBD, BwaD,
dawa, dawB, dBwa, dBwb$.

\proof
Clearly one can find such words in $P_3$.  So assume by induction
that we have $awa, awB, Bwa, BwB \in \ds_n (P_3)$.  Then 

$$ \eqalign {
[awa,awB] &= awaawBAWAbWA \cr 
&= awaawa(ABA)WAbWA \cr
&= awaawawBawB(ABA) \cr
&= awBD} \qno$$
and similarly
$$[awa,BwB] =  awaBwa(ABA)WAbWb = awaD \qno$$
$$[Bwa,awB] = Bwaawa(ABA)WbbWA = BwBD \qno$$
$$[Bwa,BwB] = BwaBwa(ABA)WbbWb =  BwaD \qno$$

The above words are all in $\ds_{n+1} (P_3)$.  Since $\ds_{n+1} (P_3)$
is a fully invariant subgroup of $P_3$, we may conjugate the
inverses of these words by $d$ to get $dawa, dawB, dBwa, dBwB$.  For example,
$(awaD)^{-1} = dAWA$, and $ddAWAD = dBwB$.  The words that end
in $D$ may now be combined with the words that begin with $d$ to
make $awa, awB, Bwa, and BwB$, and these may then be used for
the next induction step.
\square

\prop
Let $x \in B_3$, and let $n>0$.
Then there exists a $x^\prime \in B_3$, which is
a word in $a = \sigma_1$ and $B = \sigma_2^{-1}$ only,
such that $x$ is congruent to $x^\prime$ modulo
$\ds_n (P_3)$.

\proof
Write $x$ as a word in $a,b,A,B$.  Replace every occurence of
$A$ by $wDAdw = wBw$, where $wD \in \ds_n (P_3)$ is one of the
words from Proposition~3.17.  Similarly, replace $b$ by
$wDbdw = waw$.
\square

\theom
Let $K$ be a knot and let $n>0$.   Then there exist an infinite number
of prime, alternating knots in the $\ds_n$-equivalence class of $K$.

\proof
The ideas are the same as in the proof of Theorem~2 in [23].
Take a braid word in $B_k$ and insert elements in $\ds_n (P_k)$
of the forms from Proposition~3.17.
so that the braid word becomes alternating.  
Insert more such elements to make the knot prime, using Menasco's theorem [16]
that an alternating prime diagram represents a prime knot.
Then insert still more such elements to create an infinite family of
knots, using the fact that an alternating diagram 
has minimal crossing number (see Murasugi [18] or Thistlethwaite [24]).
It suffices to work
with three consecutive strands at a time, and to use the elements
of $\ds_n(P_3)$ from Proposition~3.17 (shifted over appropriately)
as in the proof of Proposition~3.22.
\square

\vfil
\eject

{\bf REFERENCES}
\bigskip
\count3 = 1

\smallskip \item{[\bnum]} D. Bar-Natan, {\it On the Vassiliev knot invariants},
Topology  {\bf 34} (1995) no. 2, 423--472

\smallskip \item {[\bnum]} D. Bar-Natan, {\it Vassiliev homotopy string link invariants},
Journal of Knot Theory and its Ramifications {\bf 4} (1995) no. 1, 13--32

\smallskip \item {[\bnum]} J. S. Birman, ``Knots, links, and mapping class groups'',
Annals of Mathematics Studies {\bf 82}, Princeton University Press, 1974

\smallskip \item {[\bnum]} J. S. Birman, {\it New points of view in knot theory},
Bulletin of the American Mathematical Society {\bf 28} (1993) no. 2, 253--287

\smallskip \item {[\bnum]} J. S. Birman and W. W. Menasco, {\it Studying links via
closed braids III: classifying links which are closed 3-braids},
Pacific Journal of Mathematics {\bf 161} (1993) no. 1, 25--113

\smallskip \item {[\bnum]} J. S. Birman and B. Wajnryb, 
{\it Markov classes in certain finite quotients of Artin's braid group}
Israel Journal of Mathematics {\bf 56} (1986) no. 2, 160--178

\smallskip \item {[\bnum]} M. Falk and R. Randell, {\it Pure braid groups and products of
free groups}, ``Braids'' (Santa Cruz, CA, 1986), 217--228,
Contemporary Mathematics {\bf 78}, American Mathematical Society, 1988

\smallskip \item {[\bnum]} S. Garoufalidis and J. Levine, {\it Finite-type $3$-manifold invariants,
the mapping class group, and blinks}, Journal of Differential Geometry {\bf 47}
(1997) no. 2, 257--320

\smallskip \item {[\bnum]} N. Gupta, ``Free group rings", Contemporary Mathematics {\bf 66},
American Mathematical Society, 1987

\smallskip \item {[\bnum]} M. N. Gusarov
{\it On $n$-equivalence of knots and invariants of finite degree},
``Topology of Manifolds and Varieties'', 173--192, Advances in
Soviet Mathematics {\bf 18}, American Mathematical Society, 1994

\smallskip \item {[\bnum]} K. Habiro, {\it Claspers and the Vassiliev skein modules},
preprint, University of Tokyo

\smallskip \item {[\bnum]} E. Kalfagianni and X.-S. Lin, {\it Regular Seifert surfaces and
Vassiliev knot invariants}, preprint {\tt GT/9804032} available from 
{\tt front.math.ucdavis.edu}

\smallskip \item {[\bnum]} T. Kohno, {\it Vassiliev invariants and the de Rham complex on the
space of knots}, ``Symplectic geometry and quantization'' (Sanda and Yokohama, 1993), 123--138,
Contemporary Mathematics {\bf 179}, American Mathematical Society, 1994

\smallskip \item {[\bnum]} X.-S. Lin, {\it Power series expansions and invariants of links}, 
``Geometric Topology'' (Athens, GA 1993), AMS/IP Studies in Advanced Mathematics {\bf 2.1},
American Mathematical Society, 1997

\smallskip \item {[\bnum]} Magnus, Karrass, and Solitar, ``Combinatorial group
theory'', Dover Publications, New York, 1976

\smallskip \item {[\bnum]} W. W. Menasco, {\it Closed incompressible surfaces in
alternating knot and link complements}, Topology {\bf 23} (1984) no. 1, 37--44

\smallskip \item {[\bnum]} J. W. Milnor, {\it Link Groups},
Annals of Mathematics {\bf 59} (1954), 177--195

\smallskip \item {[\bnum]} K. Murasugi, {\it Jones polynomials and classical
conjectures in knot theory}, Topology {\bf 26} (1987), no. 2, 187--194
 
\smallskip \item {[\bnum]} K. Y. Ng and T.~Stanford, {\it On Gusarov's groups of knots}, to appear
in the Mathematical Proceedings of the Cambridge Philosophical Society

\smallskip \item {[\bnum]} Y. Ohyama, {\it Vassiliev invariants and similarity of
knots}, Proceedings of the American Mathematical Society {\bf 123} (1995) no. 1, 287--291

\smallskip \item {[\bnum]} E. Rips, {\it On the fourth integer dimension subgroup}, 
Israel Journal of
Mathematics {\bf 12} (1972), 342--346

\smallskip \item {[\bnum]} R. Sandling, {\it The dimension subgroup problem}, Journal of Algebra
{\bf 21} (1972), 216--231

\smallskip \item {[\bnum]} T. Stanford, {\it Braid commutators and Vassiliev invariants},
Pacific Journal of Mathematics {\bf 174} (1996) no. 1, 269--276

\smallskip \item {[\bnum]} M. B. Thistlethwaite, 
{\it A spanning tree expansion of the Jones polynomial},
Topology {\bf 26} (1987), no. 3, 297--309

\smallskip \item{[\bnum]} V. A. Vassiliev, {\it Cohomology of knot spaces},
``Theory of Singularities and Its Applications", 23--69,
Advances in Soviet Mathematics {\bf 1}, American Mathematical Society, 1990

\smallskip \item {[\bnum]} B. Wajnryb, {\it Markov classes in certain 
finite symplectic representations
of braid groups}, ``Braids'' (Santa Cruz, CA, 1986), 687--695, Contemporary Mathematics
{\bf 78}, American Mathematical Society, 1988

\end